\documentclass[a4paper]{article}
\usepackage{epsfig}
\usepackage{amsmath}
\usepackage{amsfonts,amssymb}
\usepackage{amsthm}
 \newcommand{\RR}{\mathbb{R}}

\let\R=\RR

\newtheorem{theorem}{Theorem}
\newtheorem{lemma}{Lemma}

\theoremstyle{definition}
\newtheorem{remark}{Remark}

\def\pa{p_\alpha}

\begin{document}
%
\title{A note on the radial solutions for the supercritical
  H\'enon equation\thanks{Partially supported by M.I.U.R., national
    project \textit{Metodi variazionali ed equazioni differenziali non
      lineari}.}}
\author{Vivina Barutello, Simone Secchi \\
  \small{Dipartimento di Matematica ed Applicazioni, Universit\`a di
    Milano--Bicocca.} \\ \small{via R.~Cozzi 53, I-20125 Milano (Italy)} \and
     Enrico Serra \\ \small{Dipartimento di Matematica, Universit\`a
       di Milano.} \\ \small{via C.~Saldini 50, I-20133 Milano (Italy)}}
\date{}
\maketitle
\medskip


\bigskip
\bigskip
\bigskip

\begin{abstract}
We prove the existence of a positive radial solution for the H\'enon
equation with arbitrary growth. The solution is found by means of a shooting
method and turns out to be an increasing function of the radial variable.
Some numerical experiments suggest the
existence of many positive oscillating solutions.

\end{abstract}

%
\section{Introduction}

In 1982, W.-M. Ni wrote the first rigorous paper, \cite{ni82}, on an
equation introduced ten years earlier by H\'enon in \cite{he} as a
model for mass distribution in spherically symmetric clusters of stars.
This equation goes now under the name of
\textit{H\'enon equation}, and was originally coupled with
Dirichlet boundary conditions:
\begin{equation}
\label{eq:HenonD}
\begin{cases}
 -\Delta u  = |x|^\alpha u^p &\text{in } B_1, \\
 u>0, &\text{in } B_1, \\
 u= 0, &\text{on } \partial B_1,
\end{cases}
\end{equation}
where $B_1 = \left\{ x \in \R^N \mid |x|<1 \right\}$, with $N\geq 3$,
$\alpha >0$ and $p>1$.

The existence of solutions to (\ref{eq:HenonD}) for $p<\frac{N+2}{N-2}= 2^*-1$
is a standard exercise in critical point theory that can be solved
by various simple approaches.

On the other hand, Ni's main result states that \eqref{eq:HenonD} has (at least) one
solution provided $p< 2^*-1+\frac{2\alpha}{N-2}$ and thus
enlarges considerably the range of solvability beyond the classical critical threshold $p=2^*-1$.
It is also simple to prove, by means of the Pohozaev identity, that for $p \geq 2^*-1 +\frac{2\alpha}{N-2}$,
problem (\ref{eq:HenonD}) has no solution.

Thus for the Dirichlet problem (\ref{eq:HenonD}) the picture is rather sharp: setting
$\pa = 2^* + \frac{2\alpha}{N-2}$, problem (\ref{eq:HenonD}) is solvable if
and only if $p< \pa-1$. In the following we will call $\pa$ the \textit{H\'enon critical exponent}.

The key observation in Ni's work is the fact that the presence of the
weight $|x|^\alpha$, which is radial and vanishes at $x=0$, allows one
to gain compactness properties when one restricts the analysis to
radial functions.  Indeed, by means of the pointwise estimate
(\cite{ni82})
\begin{equation}
\label{grow}
|u(|x|)| \leq C \frac{\| \nabla u\|_2}{|x|^{\frac{N-2}{2}}}, \quad \hbox{for almost every $x \in B_1$},
\end{equation}
which holds true for any radial~$u\in H_0^1(B_1)$, one can easily prove that
the embedding of $H_{0,\mathrm{rad}}^1(B_1)$ into $L^p(|x|^\alpha\,dx)$ is compact precisely for
$p < \pa-1$,
and this is what one needs to prove existence of a solution to (\ref{eq:HenonD}).

Quite recently, much attention has been devoted to
\textit{symmetry--breaking} issues, namely to
the question of whether least--energy solutions of \eqref{eq:HenonD} are
radially symmetric functions (for example for large $\alpha$). In their seminal
papers \cite{ssw,SW}, Smets \textit{et al.} proved that this is indeed false
for $\alpha$ sufficiently large. In the last few years various aspects of the
H\'enon equation have been analyzed,  and the resulting literature is
nowadays rather rich (see for example \cite{BS}, \cite{bst}, \cite{BW1},
\cite{cst}, \cite{CT}, \cite{CP}, \cite{PS}, \cite{Se}, \cite{SS} and references therein).
All these papers concern the Dirichlet problem.
\medskip

The H\'enon equation just recently started
to draw attention when coupled with \textit{Neumann} boundary conditions. In this case
the problem reads
\begin{equation}\label{eq:Henon}
\begin{cases}
 -\Delta u +u = |x|^\alpha u^p &\text{in } B_1, \\
 u>0, &\text{in } B_1, \\
 \frac{\partial u}{\partial \nu}= 0, &\text{on } \partial B_1,
\end{cases}
\end{equation}
and has been studied in \cite{gs}, where the authors proved
some symme\-try--breaking results by connecting the question of symmetry
of the ground states to the symmetry properties of extremal functions in some
trace inequalities. In the paper \cite{gs} of course it is assumed that $p < 2^*-1$,
to use a variational approach in $H^1(B_1)$.

From the point of view of the mere \textit{existence} of solutions for the Neumann problem
(\ref{eq:Henon}), the situation presents both analogies and discrepancies with respect
to the Dirichlet problem. Indeed, also in the Neumann case it is very easy to check that
the problem admits at least one solution if $p < 2^* - 1$, and it has been proved
in \cite{gs} that Ni's result extends to (\ref{eq:Henon}): by using an $H^1$ version
of inequality (\ref{grow}), it is simple to prove that (\ref{eq:Henon}) admits at least
one (radial) solution for every $p< \pa-1$.

If one wishes to complete the picture of the solvability
for \eqref{eq:Henon} as a function of $p$, like in the Dirichlet case, one has to face the fact that the
Pohozaev identity gives no relevant information in presence of Neumann boundary conditions. To
our knowledge, it is not known whether the critical H\'enon exponent serves as a threshold
between existence and nonexistence in \eqref{eq:Henon}.
\medskip

The purpose of this note is to fill this gap, by looking
for radial solutions of~\eqref{eq:Henon} without any limitation on $p$.
Of course without bounds on $p$ we cannot make use of variational arguments, and we
take instead an ODE viewpoint. Our main result is the following.

\begin{theorem}
\label{th:main}
For every $p>1$ and $\alpha >0$, problem \eqref{eq:Henon} admits a
strictly increasing radial solution.
\end{theorem}

The preceding result shows that
the uselessness of the Pohozaev argument for nonexistence of solutions is
not a technical obstruction, but reflects a completely different situation
with respect to the Dirichlet problem. Though our arguments are rather simple,
we point out that it is rather difficult to find in the literature existence results
for elliptic equations
without any growth condition; some exceptions, for singularly
perturbed elliptic problems can be found
for example in \cite{amn1,amn2}.

Our main result can be easily extended to the more general problem
\begin{equation}\label{eq:general}
\begin{cases}
 -\Delta u +u = \phi(|x|) f(u), &\text{in } B_1, \\
 u>0, &\text{in } B_1, \\
 \dfrac{\partial u}{\partial \nu} = 0, &\text{on } \partial B_1,
\end{cases}
\end{equation}
under suitable assumptions on $\phi$ and $f$, but always with no growth restrictions
on $f$.

The last section of the paper contains some numerical experiments that
suggest a rather surprising fact: there are choices of $N$, $\alpha$
and $p$ for which many (and possibly infinitely many) radial solutions of \eqref{eq:Henon} seem
to exist. We clearly state that this is only a numerical hint towards further research, since no
rigorous proof has been written so far.
By the way, we believe that for any choice of the parameters, problem \eqref{eq:Henon}
has exactly one radial, increasing solution, but we have to state this only as a
\medskip

\noindent \textbf{Conjecture.} \textit{For every $N \geq 3$, $\alpha >0$
  and $p>1$, the only solution of \eqref{eq:Henon} which is positive,
  radial and increasing is that of Theorem~\ref{th:main}.}
\medskip

The matter of uniqueness for radially symmetric solutions of semilinear
elliptic equations on balls or annuli is a classical and often
overwhelmingly difficult issue. We refer to \cite{felmer} and the references therein
for a short summary of known results. We have been unable to find
uniqueness results in the literature, concerning nonautonomous
equations like \eqref{eq:Henon} with a monotonically increasing
dependence on $|x|$.

\section{An existence result for the H\'{e}non equation}
\label{sec:henon}
A radial solution for \eqref{eq:Henon}
must solve the ODE problem
\begin{equation}\label{eq:Henon_radiale}
\begin{cases}
 -u'' -\frac{N-1}{r}u'+u = r^\alpha u^p, &\text{in $(0,1)$},
 \\
 u>0, &\text{in $(0,1)$}, \\
 u'(0) = u'(1) = 0. &
\end{cases}
\end{equation}
We observe that $u'(1)=0$ corresponds to the Neumann boundary condition,
while we require $u'(0)=0$ to obtain classical solutions.
Since we do not impose any upper bound on $p$, variational techniques
do not seem to be useful to prove any existence result. For this reason we
use a shooting method, which consists in finding $\gamma>0$ such that
the solution $u_\gamma$ of the initial value problem
\begin{equation}\label{eq:shooting}
\begin{cases}
 -u'' -\frac{N-1}{r}u'+u = r^\alpha u^p, &\text{in $(0,1)$}, \\
 u>0, &\text{in $(0,1)$} \\
 u(0) = \gamma,\;  u'(0) = 0.&
\end{cases}
\end{equation}
satisfies $u'_\gamma(1)=0$. Equation \eqref{eq:shooting}$_1$ can be
written in the form
\begin{equation}\label{eq:riscritta}
 \left( r^{N-1}u'\right)' = r^{N-1}(u-r^\alpha u^p),
\end{equation}
and it is quite natural to introduce the auxiliary function $A(r) =
r^{N-1}u'(r)$.  By the definition of $A$ and \eqref{eq:shooting},
we deduce that $A(0) = 0$
and that $A$ increases strictly if and only if $u(r)<c(r)$ where
$$
c(r)=\frac{1}{r^{\alpha/(p-1)}}.
$$
The curve $c$ will play a crucial r\^{o}le
in our discussion; its peculiarity is that it is asymptotic to the
coordinate axes, i.e.
\begin{equation}\label{eq:c}
 \lim_{r \to 0+} c(r) = +\infty \quad \text{and} \quad \lim_{r \to +\infty} c(r) = 0^+.
\end{equation}
The next Lemma is well known, see for example \cite{caku}.
\begin{lemma} \label{lem:1}
  For every $\gamma >0$, problem \eqref{eq:shooting} is uniquely
  solvable on $[0,+\infty)$. Its solution $u_\gamma$ is continuously differentiable
  with respect to the initial value $\gamma$.
\end{lemma}

We begin now a qualitative study of solutions to \eqref{eq:shooting},
with the aim of  proving the existence of (at least) an initial datum
$\bar\gamma >0$ such that the corresponding $u_{\bar \gamma}$ matches
the Neumann boundary condition at $r=1$.

Take any~$\gamma >0$, and consider the solution $u_\gamma$ of
\eqref{eq:shooting}.  The corresponding auxiliary function
$A_\gamma(r) = r^{N-1}u'_\gamma(r)$ is strictly positive from $r=0$
until $u_\gamma$ intersects the curve $c$,
see~\eqref{eq:riscritta}. After the first crossing, the derivative
of $A_\gamma$ becomes negative. Therefore there exists a zero of $A_\gamma'$, and the
first intersection point
\[
 r_\gamma := \inf \{r >0 \mid A'_\gamma (r) = 0 \}
           = \inf \{r >0 \mid u_\gamma (r) = c(r) \}
\]
between the graph of $u_\gamma$ and that of $c$ is well defined.  By
definition, $u'_\gamma (r)>0$ (at least) on the interval
$(0,r_\gamma)$.  Let us call $R_\gamma$ the first nontrivial stationary point of
$u_\gamma$, namely
\[
 R_\gamma := \inf \{r >0 \mid u'_\gamma(r)=0 \}.
\]
The next Lemma states in particular that $R_\gamma$ is well defined.

\begin{lemma}\label{lem:2}
For every $\gamma > 0$ we have $r_\gamma < R_\gamma < +\infty$.
\end{lemma}
\begin{proof}
Since the function $A_\gamma$ steadily increases on $(0,r_\gamma)$
and $A(0)=0$, we deduce that $A_\gamma(r_\gamma)>0$ and then that
$u'_\gamma(r_\gamma)>0$. Hence $r_\gamma < R_\gamma$.  Let us show
that $R_\gamma < + \infty$. Assume, for the sake of contradiction,
that $R_\gamma = + \infty$. This means that $u'(r)>0$ for all $r
\geq r_\gamma$. Fix any $R_0 > r_\gamma$ (of course $u'_\gamma
(R_0)>0$) and choose $\delta$ such that
\[
0< \delta < \left( \frac{u_\gamma(R_0)}{c(R_0)} \right)^{p-1}-1.
\]
This choice is possible since $u_\gamma$ increases strictly
and $c$ decreases strictly, so that $u_\gamma(R_0)>c(R_0)$.
Therefore
\begin{equation} \label{eq:7} r^\alpha u\sb\gamma^{p-1}(r)-1 >
  \delta >0 \quad\hbox{for every } r \geq R_0.
\end{equation}
We now integrate \eqref{eq:shooting}$_1$ on the interval $[R_0,R]$, with
$R>R_0$; we have
\begin{eqnarray*}
0 &=& u'_\gamma(R)-u'_\gamma(R_0) + (N-1) \int_{R_0}^R \frac{u'_\gamma (r)}{r}\,dr\\
&&\qquad \quad- \int_{R_0}^R u_\gamma(r)\,dr + \int_{R_0}^R r^\alpha u_\gamma^p(r)\,dr\\
&\geq& -u'_\gamma(R_0)-\int_{R_0}^R u_\gamma(r) \left[1-r^\alpha u_\gamma^{p-1}(r)\right] \,dr
\end{eqnarray*}
and hence
\[
\int_{R_0}^R u_\gamma (r) \left[ r^\alpha u_\gamma (r)^{p-1} - 1 \right]\,dr
\leq u'_\gamma(R_0)
< +\infty.
\]
for every $R >R_0$. However, thanks to \eqref{eq:7}, we obtain
\[
\lim_{R \to +\infty} \int_{R_0}^R u_\gamma (r) \left[ r^\alpha u_\gamma (r)^{p-1} -1 \right] \, dr
\geq \lim_{R \to +\infty} (R-R_0) \gamma \delta = +\infty,
\]
a contradiction that concludes the proof.
\end{proof}

\begin{lemma} \label{lem:3}
There results $u''_\gamma (R_\gamma) < 0$.
\end{lemma}
\begin{proof}
Indeed, from equation \eqref{eq:shooting} we obtain
\[
-u''_\gamma (R_\gamma) = u_\gamma (R_\gamma)
\left( 1-R_\gamma^\alpha u^{p-1}_\gamma (R_\gamma) \right)
\]
which is a positive quantity since $u_\gamma (R_\gamma) > c(R_\gamma)$.
\end{proof}

\begin{lemma} \label{lem:4}
  There exists $\delta>0$ such that for every $\gamma<\delta$ there
  results~$R_\gamma>1$.
\end{lemma}
\begin{proof}
  Lemma \ref{lem:1} implies that $\sup_{r \in [0,1]} |u_\gamma (r)| <
  1$ for all $\gamma$ sufficiently small. Therefore $u_\gamma$ lies
  below the curve $c$ on $[0,1]$, since $c$ decreases and $c(1)=1$
  independently of the parameters $p$ and $\alpha$. The claim follows
  from Lemma \ref{lem:2}.
\end{proof}

\begin{lemma} \label{lem:5}
  $\lim_{\gamma \to +\infty}R_\gamma = 0$.
\end{lemma}
\begin{proof}
  For the sake of contradiction we suppose that there exist $\delta >
  0$ and a sequence $(\gamma_k)_k$, $\gamma_k \to +\infty$, such that
  $R_{\gamma_k} \geq \delta$ for every $k$.  Since
  $u_{\gamma_k}'$ is strictly positive on $(0,\delta/2)$, the function
  $A_{\gamma_k}$ is strictly positive on $(0,\delta/2]$.  Hence, for
  every $k$,
\[
\begin{split}
 0 < A_{\gamma_k}(\delta/2) & = \int_0^{\delta/2} A'_{\gamma_k}(r)\,dr \\
 & = \int_0^{\delta/2} r^{N-1}u_{\gamma_k}(r)\left(1-r^\alpha u_{\gamma_k}^{p-1}(r)\right)\,dr \\
 & \leq \int_0^{\delta/2} r^{N-1}u_{\gamma_k}(r) \left(1-r^\alpha
    \gamma_k^{p-1}\right)\,dr,
\end{split}
\]
where the last inequality holds because $u_{\gamma_k}(r) \geq \gamma_k$ on $[0,\delta/2]$.
Since $\gamma_k \to +\infty$, we can choose $k_0$ such that, for every $k \geq k_0$,
there results $\gamma_k^{-(p-1)/\alpha}<\delta/2$; for such values of $k$ we can split the integral
to obtain
\[
\begin{split}
0 &< \int_0^{\gamma_k^{-(p-1)/\alpha}} r^{N-1}u_{\gamma_k}(r)
       \left(1-r^\alpha \gamma_k^{p-1}\right)\,dr\\
& \qquad\qquad+ \int_{\gamma_k^{-(p-1)/\alpha}}^{\delta/2} r^{N-1}u_{\gamma_k}(r)
       \left(1-r^\alpha \gamma_k^{p-1}\right)\,dr\\
& \leq u_{\gamma_k}(\gamma_k^{-(p-1)/\alpha})\int_0^{\gamma_k^{-(p-1)/\alpha}} r^{N-1}
       \left(1-r^\alpha \gamma_k^{p-1}\right)\,dr\\
&\qquad\qquad + u_{\gamma_k}(\gamma_k^{-(p-1)/\alpha})\int_{\gamma_k^{-(p-1)/\alpha}}^{\delta/2} r^{N-1}
       \left(1-r^\alpha \gamma_k^{p-1}\right)\,dr\\
&= u_{\gamma_k}(\gamma_k^{-(p-1)/\alpha}) \int_0^{\delta/2} r^{N-1}
       \left(1-r^\alpha \gamma_k^{p-1}\right)\,dr.
\end{split}
\]
Indeed $u_{\gamma_k}$ increases and the quantity $(1-r^\alpha
\gamma_k^{p-1})$ is positive on the interval $(0,\gamma_k^{-(p-1)/\alpha})$ and
negative on the interval $(\gamma_k^{-(p-1)/\alpha},\delta/2)$. Integrating we reach
the contradiction
\[
0 < u(\gamma_k^{-(p-1)/\alpha}) \left(\frac{\delta}{2}\right)^{N}
\left[\frac{1}{N}-\frac{\gamma_k^{p-1}}{N+\alpha}\left(\frac{\delta}{2}\right)^\alpha
\right] < 0
\]
whenever $\gamma_k >\left( \frac{N+\alpha}{N}
\right)^{1/(p-1)}(\delta/2)^{-\alpha/(p-1)}$, namely for $\gamma_k$ large enough.
\end{proof}

We are now ready to state our main result.

\setcounter{theorem}{0}
\begin{theorem} \label{main}
For every $p>1$ and $\alpha >0$, problem \eqref{eq:Henon} admits a
strictly increasing radial solution.
\end{theorem}
\begin{proof}
  We have seen that the first critical point $R_\gamma$ of $u_\gamma$
  is larger than 1 for small values of $\gamma$, and smaller than 1
  for large values of $\gamma$.  Consider the map $F \colon
  (0,+\infty) \times (0,+\infty) \to \mathbb{R}$ defined by
  $F(\gamma,r)=u'_\gamma(r)$. We have that $F(\gamma,R_\gamma)=0$,
  and Lemma \ref{lem:3} implies that
\[
\partial_2 F (\gamma,R_\gamma) = u''_\gamma (R_\gamma) < 0.
\]
The Implicit Function Theorem shows that $\gamma \mapsto R_\gamma$ is
a continuous and even differentiable function. Therefore, there exists
$\bar\gamma$ such that $R_{\bar \gamma}=1$. This means that
$u_{\bar\gamma}$ is a radial solution of the Neumann problem
\eqref{eq:Henon}.
\end{proof}

\begin{remark}
When $p$ satisfies the condition
\begin{equation}\label{eq:p-1}
p-1 < \frac{\sqrt{(N-2)^2+4}-(N-2)}{2}\alpha,
\end{equation}
the Maximum Principle shows that every radial solution for the problem
\eqref{eq:Henon} can change its monotonicity at most once,
since it can intersect the curve $c$ at most twice. Indeed, computing
\begin{multline*}
-c''(r)-\frac{N-1}{r}c'(r)+c(r) =\\
(p-1)^2 r^{-2-\alpha/(p-1)}\left[r^\alpha + \frac{\alpha(N-1)(p-1) -
\alpha(\alpha+p-1)}{(p-1)^2}\right],
\end{multline*}
it is easy to see that since $r \in [0,1]$, condition \eqref{eq:p-1} implies that $c$
is a subsolution for the operator $-\Delta + I$.
Furthermore, a positive solution $u$ for \eqref{eq:Henon} is a supersolution
for the same operator. Therefore $w=u-c$ is a supersolution that vanishes when
$u$ intersects $c$. From the Maximum Principle we deduce that $w$ can vanish
at most twice.
\end{remark}

\section{A generalization}
\label{sec:general}
Theorem \ref{main} can be adapted to problem \eqref{eq:general}
by imposing some suitable assumptions of the functions $\phi$ and $f$.
First of all, both functions $\phi$ and $f$ are defined (and continuous)
on $[0,+\infty)$,
since we are looking for positive radial solutions;
furthermore we require:
\begin{description}
 \item[(h1)] $\phi$ is increasing, $\phi(0)=\ell \geq 0$ and
$\displaystyle \lim_{r \to +\infty} \phi(r)=\kappa \in [l,+\infty]$;
 \item[(h2)] the function $s \mapsto f(s)/s$ is strictly increasing;
 \item[(h3)] $\displaystyle \lim_{s \to +\infty} \dfrac{f(s)}{s} =
   \dfrac{1}{\ell}$ ($=+\infty$ if $\ell=0$);
 \item[(h4)] $\displaystyle \lim_{s \to 0^+} \dfrac{f(s)}{s} =
   \dfrac{1}{\kappa}$ ($=0$ if $\kappa=+\infty$).
\end{description}
Under conditions {\bf (h1)}--{\bf (h4)} the equation
${u}/{f(u)}= \phi(r)$ defines implicitly a continuous curve $u=\xi(r)$
which plays the same r\^{o}le as the curve $c$ in the previous
section.  Indeed if we call $H(u)={u}/{f(u)}$, then
$\xi(r)=H^{-1}\left( \phi(r) \right)$, which decreases since $H^{-1}$
decreases and $\phi$ increases.  Furthermore $\xi$ is asymptotic to
the coordinate axes. Indeed from {\bf (h3)} we get
\[
 \lim_{r\to 0+} \xi(r) = \lim_{r\to 0+}H^{-1}\left( \phi(r) \right)
                       = \lim_{u\to l}H^{-1}(u) = +\infty;
\]
similarly, from {\bf (h4)},
\[
 \lim_{r\to +\infty} \xi(r) = \lim_{r\to +\infty}H^{-1}\left( \phi(r) \right)
                       = \lim_{u\to k}H^{-1}(u) = 0.
\]
We can now proceed exactly as in Section \ref{sec:henon} defining
the auxiliary function $A$ and the points $r_\gamma$, $R_\gamma$.
Lemmas \ref{lem:2}, \ref{lem:3}, \ref{lem:4} and \ref{lem:5}
can be proved also in this setting with  minor changes. We then obtain
the required generalization:
\begin{theorem} \label{main2}
Let $\phi$, $f:[0,+\infty)\to\RR$ be continuous functions.
If assumptions {\bf (h1)}--{\bf (h4)} are satisfied then
problem \eqref{eq:general} admits a strictly increasing radial solution.
\end{theorem}

\begin{remark}\label{esempi}
  Our assumptions are clearly satisfied by nonlinearities with arbitrarily fast
  growth at infinity, like $f(s)=\exp (s)-1$, or $f(s) = \exp (\gamma s^q)-1$ for
  $\gamma>0$ and $q>1$. This latter case is particularly interesting, because it
  corresponds to Trudinger--Moser type problems without any restriction
  on $q$ and $\gamma$. Though in this paper we work in dimension
  $N\geq 3$, it is immediate to check that our results hold also for
  $N=2$, which is the case of the Trudinger--Moser problem.

  Of course,
  \textbf{(h3)} requires~$\ell=0$, i.e. $\phi (0)=0$. Similarly, a
  nonlinearity that is superlinear at zero forces $\kappa$ to be
  infinite. Also, we notice that homogeneity of $\phi$ plays no role, as long as
  $\phi $ satisfies the above assumptions; for example $\phi(r) = r^\alpha + r^\beta$,
  with $\alpha$, $\beta>0$,
  is an admissible function.
\end{remark}

\begin{remark}
It is proved in \cite{adya91} that in dimension $N\geq 3$, for any $p>1$
and for large values of $R$, the problem
\[
\begin{cases}
-\Delta u + u = u^p &\hbox{for $|x|<R$}\\
\frac{\partial u}{\partial \nu} =0 &\hbox{for $|x|=R$}
\end{cases}
\]
does not have any positive radial solution whose derivative changes sign.
We could not find any similar statement for a nonlinearity
like $|x|^\alpha u^p$.
\end{remark}

\section{Some numerical results}
This section is devoted to the description of some numerical experiments.
Such results are purely numerical and non rigorous; the purpose of the
authors is, on one hand, to give some examples of the existence result
proved in Section \ref{sec:henon}. On the other hand we want to point out
some features of interest in the behavior of the solutions for the shooting problem
\eqref{eq:shooting} when $\gamma$ diverges to $+\infty$. These numerical experiments seem to
indicate that the structure of the set of radial solutions of problem (\ref{eq:Henon})
is still far to be understood, and deserves further study.

\subsection{The monotone solution}

The monotone solution for the H\'enon equation corresponds to
a choice of the parameter $\gamma$ such that the first maximum point
$R_\gamma$ of the solution of the shooting equation
coincides with $1$.

\noindent \begin{minipage}{0.5\textwidth}
In the table on the right we have collected
some values of $\gamma$, depending on $N$, $\alpha$ and $p$,
for which $|R_\gamma - 1|<10^{-6}$.
The starred values of $p$ are the H\'enon critical exponents,
that is $p= \pa-1$.
In the first four rows of the table, we fix the dimension $N$
and the exponent $\alpha$
and we choose three values of $p$:
subcritical, critical and supercri\-ti\-cal.
We observe that $\gamma$ seems to be a decreasing function of $p$.

In the last row we have investigated the behavior of $\gamma$ as $\alpha$
becomes larger and larger while $p$ has a supercritical value.
The numerical results seem to show that $\gamma$ continues to lie near $1$.
This is probably due to the fast convergence of $R_\gamma$ to $0$ as
$\gamma \to +\infty$, see Lemma~\ref{lem:5}.
\end{minipage}
\hfill
\begin{minipage}{0.45\textwidth}
\[
\begin{array}{|c|c|c|c|}
\hline
N & \alpha & p  & \gamma \\[3pt]
\hline\hline
           &    &  5   & 1.0816 \\[3pt]
\;3\;      &  3 & 11^* & 0.9710 \\[3pt]
           &    & 15   & 0.9487 \\[3pt]
\hline
  &        &  3   & 1.3739 \\[3pt]
4 &   5    &  8^* & 1.0306 \\[3pt]
  &        & 12   & 0.9872 \\[3pt]
\hline
  &        &  4         & 1.3102 \\[3pt]
5 &   9    & \;25/3^*\; & 1.0632 \\[3pt]
  &        & 12         & 1.0147 \\[3pt]
\hline
   &       &  11/4    & 1.2175 \\[3pt]
10 &   5   &  5^*     & 1.0688 \\[3pt]
   &       &  10      & 1.0105 \\[3pt]
\hline\hline
         &   50    &  20  & 1.0485    \\[3pt]
\; 10 \; &  100    &  50  & 1.0114    \\[3pt]
         &\;200 \; &  50  &\;1.0135\; \\[3pt]
\hline
\end{array}
\]
\end{minipage}

\subsection{Numerical evidence of oscillating solutions}
We have considered so far the existence of a radial solution
with Neumann boundary conditions with the \emph{first} stationary
point at $1$.
Recalling Ni's results about the existence of oscillating
radial solutions for some elliptic problem on $\RR^n$ (see \cite{ni83})
we now investigate the existence on non--monotone solutions for our problem.
We address two natural questions. Can the shooting solution,
which is defined on the whole interval $[0,+\infty)$,
have stationary points different from $R_\gamma$? Can we choose $\alpha$, $p$,
and $\gamma$ such that one of such points coincides with $1$?
Although we do not have any rigorous proof of these facts,
some numerical experiments show that
the answers to these questions strongly depend on the parameters $\alpha$ and $p$,
so that it seems very unlikely to obtain a single general result.

\begin{figure}
\begin{center}
  {\psfig{figure=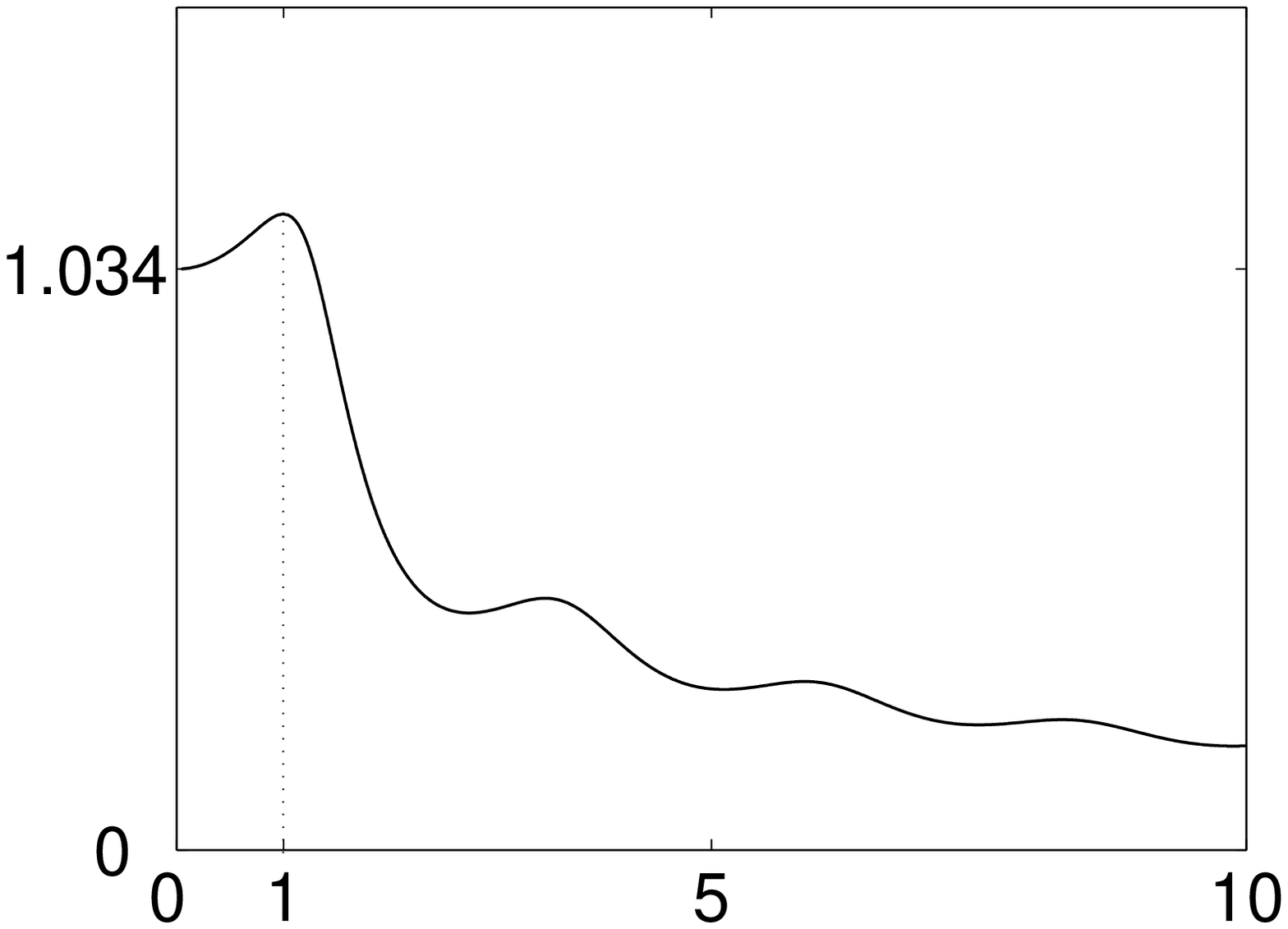,width=6.0cm}}
  {\psfig{figure=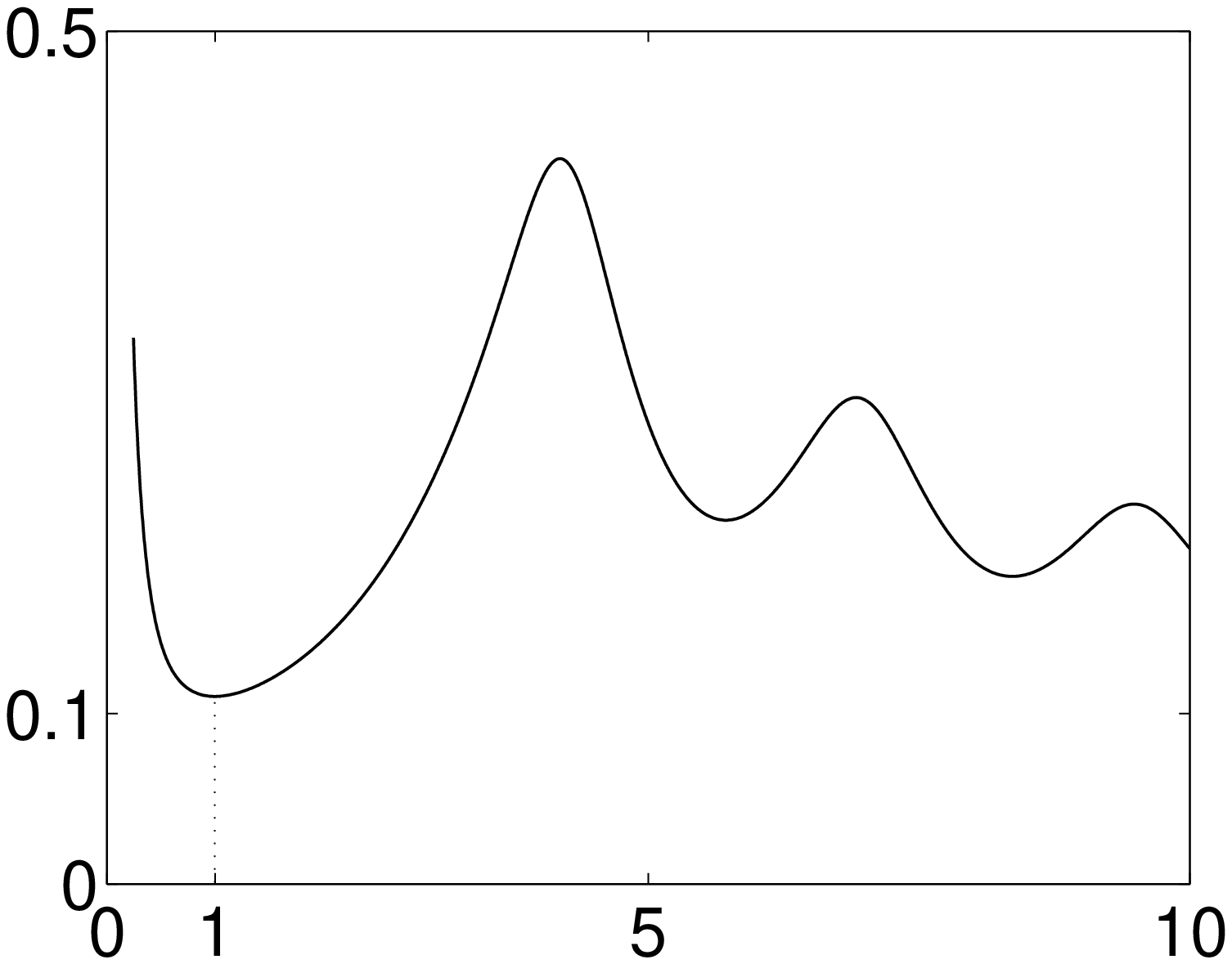,width=6.0cm}}
  {\psfig{figure=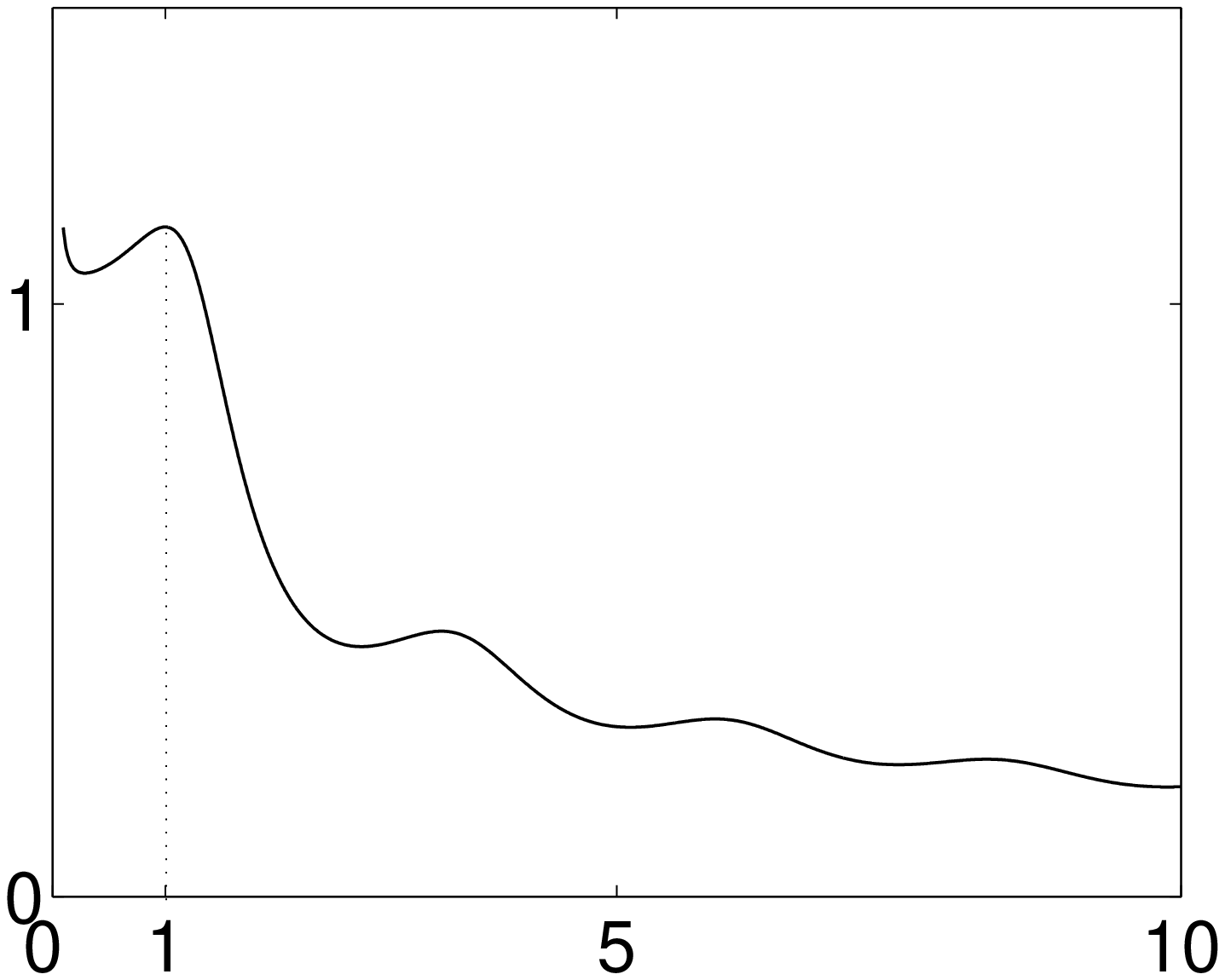,width=6.0cm}}
\end{center}
\caption{the three plots represent the shooting solution to problem \eqref{eq:shooting},
corresponding to different values of $\gamma$,
when $N=4$, $\alpha=5$ and $p= 8$ (critical). In the horizontal axis the radial variable $r$
varies to the interval $[0,10]$.
In the figure above--left we choose $\gamma=1.034$ and we obtain the first
stationary point $R_\gamma \approx 1$; above--right, when $\gamma \approx 155$,
the second stationary point approximates $1$. If $\gamma \approx 2584$
the third stationary point satisfies $R^3_\gamma \approx 1$. In the second and third picture,
the first maximum point $R_\gamma$ is not printed in order to obtain
a reasonable scaling, and we just plot the solution for $r>0.2$.
}\label{figura1}
\end{figure}

\begin{figure}
\begin{center}
  {\psfig{figure=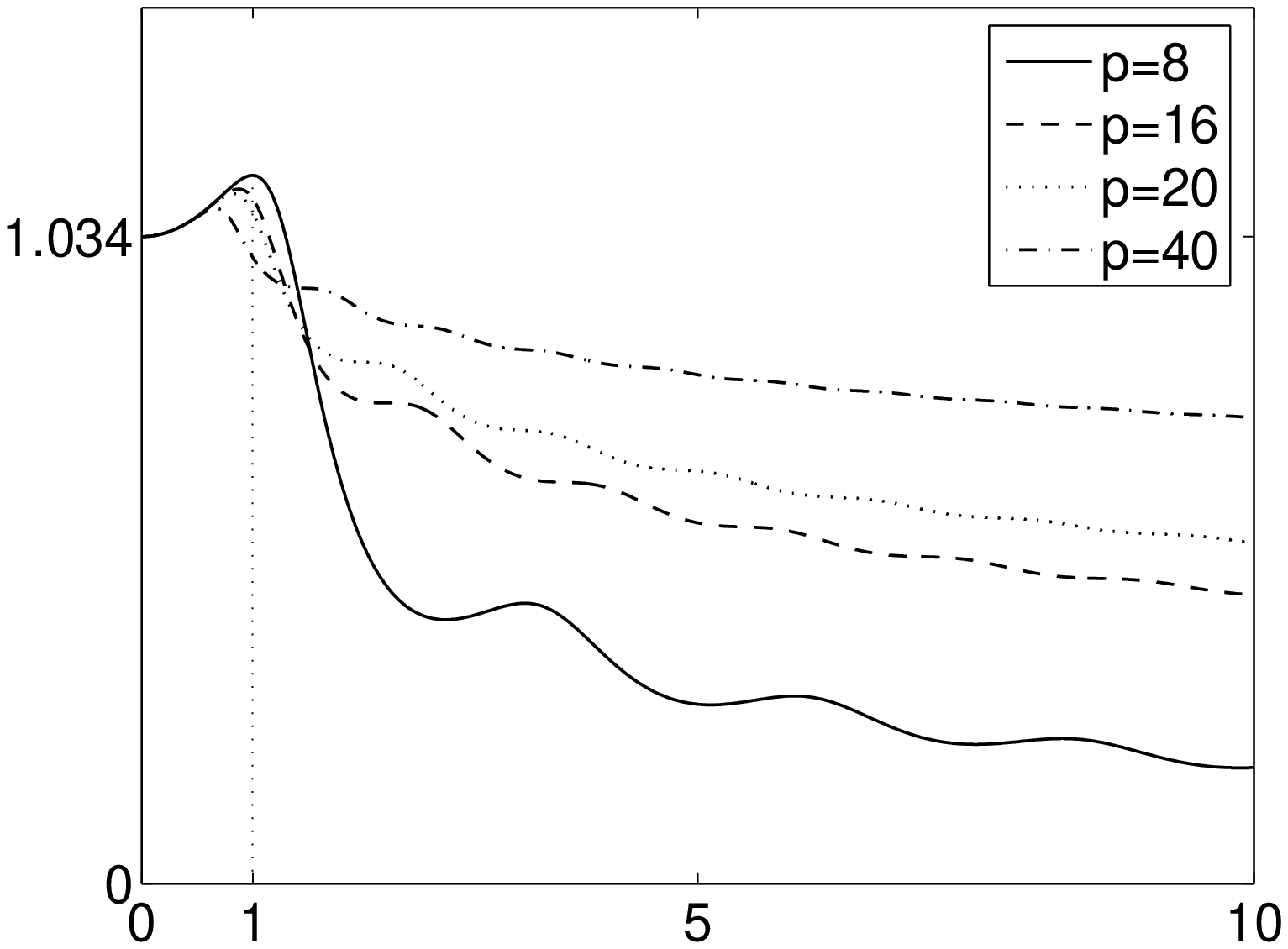,width=10.0cm}}
\end{center}
\caption{in the picture we compare the plots of the shooting solutions
when $N=4$, $\alpha=5$ and $\gamma=1.034$, when $p$ varies from the critical H\'enon
exponent $p=8$ to some supercritical values.
In the horizontal axis the radial variable $r$
varies in the interval $[0,10]$. When $p$ is critical the numerical
solution oscillates sharply; as $p$ increases the oscillations become weaker
and weaker and they disappear when $p>16$.
}\label{figura2}
\end{figure}

Fix for instance $N=4$ and $\alpha=5$. The critical H\'enon exponent $\pa-1$
is in this case
$p=8$. For such values of the parameters, the monotone solution corresponds
to $\gamma\approx1.034$. If we compute the solution on a larger interval we
observe that it oscillates. Let us call $R^n_\gamma$ the $n$--th
stationary point of this solution. Lemma \ref{lem:5} suggest that,
as $\gamma \to +\infty$, $R^n_\gamma$ decreases also for $n>1$;
for this reason we increase $\gamma$ in order to obtain $R_\gamma^n=1$
for some $n>1$.
To obtain $R^2_\gamma \approx 1$ we need to choose
$\gamma \approx 155$, while for $R^3_\gamma \approx 1$, a value of
$\gamma \approx 2584$ will do (see Figure \ref{figura1}). These cases point towards the possibility
of the coexistence of multiple radial positive solutions.
Rather surprisingly this interesting behavior gently disappears when $p$ becomes
supercritical. Indeed the oscillations become less and less sharp as $p$ increases
and disappear when $p$ is greater than some $\bar p$ (in the described case
$\bar p \approx 16$). This behavior is illustrated in Figure \ref{figura2}
and it suggests a strong difference with the autonomous case studied by Ni
in \cite{ni83}.

\end{document}